\documentclass[12pt,a4paper]{article}
\usepackage[cp1251]{inputenc}
\usepackage[T2A]{fontenc}
\usepackage[english]{babel}

\usepackage{amsmath,amsfonts,amssymb,mathrsfs,amsopn,indentfirst,amscd}

\usepackage{graphicx}
\usepackage{amsthm}
\usepackage{hyperref}

\usepackage{amssymb}
\usepackage{amsmath}
\newcommand{\ver}{{\rm ver}}
\newcommand{\vo}{{\rm vol}}

\newtheorem*{corollary*}{Corollary}

\begin{document}

\title{
On Some Problems \\
Related to a Simplex and a Ball}
\author{Mikhail Nevskii\footnote{Department of Mathematics,
              P.G.~Demidov Yaroslavl State University, Sovetskaya str., 14, Yaroslavl, 150003, Russia 
              orcid.org/0000-0002-6392-7618 
              mnevsk55@yandex.ru} 
}       
\date{May 5, 2019}
\maketitle

\begin{abstract}
Let $C$ be a convex body and let $S$ be a nondegenerate
simplex in ${\mathbb R}^n$. 
Denote by $\xi(C;S)$  the minimal 
$\tau>0$ such that $C$ is a subset of the simplex $\tau S$.
By $\alpha(C;S)$ we mean the minimal $\tau>0$ such that
$C$ is contained in a translate of $\tau S$. Earlier the author
has proved the equalities
$\xi(C;S)=(n+1)\max\limits_{1\leq j\leq n+1}
\max\limits_{x\in C}(-\lambda_j(x))+1$  \ (if $C\not\subset S$), \ 
$\alpha(C;S)=
\sum\limits_{j=1}^{n+1} \max\limits_{x\in C} (-\lambda_j(x))+1.$
Here $\lambda_j$ are linear functions called the
basic Lagrange polynomials corresponding to $S$.
In his previous papers, the author has investigated  
these formulae if
$C=[0,1]^n$. 
The present paper is related to the case
when
$C$ coincides with the unit Euclidean ball $B_n=\{x: \|x\|\leq 1\},$ where
$\|x\|=\left(\sum\limits_{i=1}^n x_i^2 \right)^{1/2}.$ We establish
various relations for 
$\xi(B_n;S)$  and $\alpha(B_n;S)$, as well as we give their
geometric interpretation.

\medskip 

\end{abstract}


\section{Preliminaries}\label{nev_s1}
 Everywhere further $n\in{\mathbb N}.$ An element
$x\in{\mathbb R}^n$ is written in the form
$x=(x_1,\ldots,x_n).$ 
By definition, 
$$\|x\|=\sqrt{(x,x)}=\left(\sum\limits_{i=1}^n x_i^2\right)^{1/2},$$ 
$$B\left(x^{(0)};\varrho\right):=\{x\in{\mathbb R}^n: \|x-x^{(0)}\|\leq \varrho \} 
\quad \left(x^{(0)}\in {\mathbb R}^n,
\varrho>0\right),$$ \ 
$$B_n:=B(0;1), \quad 
Q_n:=[0,1]^n, \quad
Q_n^\prime:=[-1,1]^n.$$

Let $C$ be a convex body in ${\mathbb R}^n$.  
Denote by $\tau C$ the image of $C$ under the homothety with center of homothety
in the center of gravity of $C$ and ratio of homothety
$\tau.$ 
For an~$n$-dimensional nondegenerate simplex $S$,
consider the value 
$\xi(C;S):=\min \{\sigma\geq 1: C\subset \sigma S\}.$
We call this number the {\it absorption index 
of $S$ with respect to $C$.}
Define $\alpha(C;S)$ as minimal $\tau>0$ such that convex body
$C$ is a subset of the simplex $\tau S$. 
By $\ver(G)$ we mean the set of vertices
of convex polytope
$G$.

Let
$x^{(j)}=\left(x_1^{(j)},\ldots,x_n^{(j)}\right),$ 
$1\leq j\leq n+1,$ be the vertices of simplex $S$.  
The matrix 
$${\bf A} :=
\left( \begin{array}{cccc}
x_1^{(1)}&\ldots&x_n^{(1)}&1\\
x_1^{(2)}&\ldots&x_n^{(2)}&1\\
\vdots&\vdots&\vdots&\vdots\\
x_1^{(n+1)}&\ldots&x_n^{(n+1)}&1\\
\end{array}
\right)$$
is nondegenerate. 
By definition, put 
${\bf A}^{-1}$ $=(l_{ij})$. 
Linear polynomials
$\lambda_j(x)=
l_{1j}x_1+\ldots+
l_{nj}x_n+l_{n+1,j}$
whose coefficients make up the columns of
${\bf A}^{-1}$
have the property  
$\lambda_j\left(x^{(k)}\right)$ $=$ 
$\delta_j^k$, where $\delta_j^k$ is the Kronecker $\delta$-symbol.
We call $\lambda_j$ the
{\it basic Lagrange polynomials corresponding to
$S$.}
The numbers $\lambda_j(x)$  
are barycentric coordinates of a point
$x\in{\mathbb R}^n$ with respect to
$S$. Simpex $S$ is given by the system of linear inequalities
$\lambda_j(x)\geq 0$. For more details about $\lambda_j$, 
see [3; Chapter\,1].

The equality $\xi(C;S)=1$ is equivalent to the inclusion
$C\subset S.$ If $C\not\subset S$, then  
\begin{equation}\label{ksi_cs_equality}
\xi(C;S)=(n+1)\max_{1\leq j\leq n+1}
\max_{x\in C}(-\lambda_j(x))+1. 
\end{equation}
(the proof was given in [2]; see also [3;\S\,1.3]). 
The relation 
\begin{equation}\label{relation_cs}
\max\limits_{x\in C} \left(-\lambda_1(x)\right)=
\ldots=
\max\limits_{x\in C} \left(-\lambda_{n+1}(x)\right)
\end{equation}
holds true if and only if 
the simplex $\xi(C;S)S$ 
is circumscribed around
convex body $C.$
In the case $C=Q_n$  equality (\ref{ksi_cs_equality}) can be reduced to the form  
$$\xi(Q_n;S)=(n+1)\max_{1\leq j\leq n+1}
\max_{x\in \ver(Q_n)}(-\lambda_j(x))+1 $$
and    (\ref{relation_cs}) is equivalent to the relation
\begin{equation}\label{relation_qs}
\max\limits_{x\in \ver(Q_n)} \left(-\lambda_1(x)\right)=
\ldots=
\max\limits_{x\in \ver(Q_n)} \left(-\lambda_{n+1}(x)\right).
\end{equation} 

For any $C$ and $S$, we have $\xi(C;S)\geq\alpha(C;S)$. The equality
$\xi(C;S)=\alpha(C;S)$ holds only in the case when
the simplex
$\xi(C;S)S$ is circumscribed around
convex body $C.$  
This is equivalent to 
(\ref{relation_cs}) and also  to 
 (\ref{relation_qs}) when $C=Q_n$.

It was proved in [4] (see also
 [3; \S\,1.4])
that
\begin{equation}\label{alpha_cs_equality}
\alpha(C;S)=
\sum_{j=1}^{n+1} \max_{x\in C} (-\lambda_j(x))+1.
\end{equation}
If $C=Q_n$, then this formula 
can be written in rather more geometric way:
\begin{equation}\label{alpha_d_i_formula}
\alpha(Q_n;S)
=\sum_{i=1}^n\frac{1}{d_i(S)}. 
\end{equation}
Here $d_i(S)$ is {\it the $i$th axial diameter of simplex $S$,} i.\,e.,
the length of a longest segment in $S$
parallel to the $i$th
coordinate axis. Equality (\ref{alpha_d_i_formula})
was obtained in [11].
When $S\subset Q_n,$ we have $d_i(S)\leq 1.$ Therefore, for these simplices,
(\ref{alpha_d_i_formula}) gives
\begin{equation}\label{ksi_alpha_n_ineq}
\xi(Q_n;S)\geq\alpha(Q_n;S)
=\sum_{i=1}^n\frac{1}{d_i(S)}\geq n.
\end{equation}

Earlier the author established the equality
\begin{equation}\label{d_i_l_ij_formula}
\frac{1}{d_i(S)}=\frac{1}{2}\sum_{j=1}^{n+1} |l_{ij}|
\end{equation}
(see [2]).
Being combined together, (\ref{alpha_d_i_formula}) and (\ref{d_i_l_ij_formula}) yield
\begin{equation}\label{alpha_qs_formula}
\alpha(Q_n;S)=\frac{1}{2}\sum_{i=1}^n\sum_{j=1}^{n+1} |l_{ij}|. 
\end{equation}
Note that $\alpha(C;S)$ is invariant under parallel translation of the sets
and for $\tau>0$ we have
$\alpha(\tau C;S)=\tau\alpha(C;S).$  Since
$Q_n^\prime=[-1,1]^n$  
is a translate of the cube $2Q_n$, 
after replacing $Q_n$  with
$Q_n^\prime$ we obtain from (\ref{alpha_qs_formula})
an even simpler formula: 
\begin{equation}\label{alpha_q_prime_s_formula}
\alpha(Q_n^\prime;S)=\sum_{i=1}^n\sum_{j=1}^{n+1} |l_{ij}|. 
\end{equation}

Let us define the value
$$\xi_n:=\min \{ \xi(Q_n;S): \,
S \mbox{ --- $n$-мерный симплекс,} \,
S\subset Q_n, \, \vo(S)\ne 0\}.$$
Various estimates of $\xi_n$ were obtained first by the author and then
by the author and A.\,Yu.~Ukhalov
(e.\,g., see papers [1], [2], [5], [6], [7], [8], [12]
and book
[3]).
Always $n\leq \xi_n<n+1$. 
Nowaday the precise values of $\xi_n$ 
are known for $n=2,5,9$ and also for the infinite set of odd $n$'s
for any of which there exists an Hadamard matrix of order $n+1$. 
If $n\ne 2$, then every known value of
$\xi_n$ is equal to $n$, whereas $\xi_2=1+\frac{3\sqrt{5}}{5}=2.34\ldots$
Still remains unknown 
is there exist
an even $n$ with the property
$\xi_n=n$. 
There are some other open problems concerning the numbers
$\xi_n$.

In this article, we will discuss 
the analogues of the above characteristics 
for a simplex and an Euclidean ball. 
Replacing a cube with a ball makes many 
questions much more simpler. However, geometric interpretation of general results
has a certain interest also in this particular case.
Besides, we will note some new applications of the basic Lagrange polynomials.

Numerical characteristics connecting
simplices and subsets of ${\mathbb R}^n$
have applications for obtaining various estimates 
in polynomial interpolation of functions  defined
on mul\-ti\-dimen\-sional domains.  This approach and the corresponding analytic
methods
in detailes
were described in 
[3]. Lately these questions have been managed to study
also by computer methods (see, e.\,g., 
[5], [6], [8], [12]).

\section{The value $\alpha(B_n;S)$}\label{nev_s2}

The {\it inradius of an $n$-dimensional simplex $S$} is the
maximum of the radii of balls contained within $S$.
The center of this unique maximum ball is called the {\it incenter of $S$.}
The boundary of the maximum ball is a sphere that has a single common point
with each $(n-1)$-dimensional face of $S$. By the {\it circumradius of S} 
we mean the minimum of the radii of balls containing $S$.
The boundary of this unique minimal ball does not necessarily contain all the
vertices of $S$. Namely, this is only when the center of the minimal ball
lies inside the simplex.

The inradius $r$ and the circumradius $R$ of a simplex $S$ 
satisfy the so-called {\it Euler inequality}
\begin{equation}\label{euler_ineq}
R\geq nr.
\end{equation}
Equality in
(\ref{euler_ineq}) takes place if and only if 
$S$ is a regular simplex.
Concerning the proofs of the  Euler inequality, its history and generalizations, 
see, e.\,g., [10], [13], [14].

In connection with 
(\ref{euler_ineq}), let us remark an analogue to the following property being
true for parallelotopes 
(see [11], [3; \S\,1.8]).
{\it Let   $S$ be a nondegenerate simplex and let
$D,$ $D^*$~be parallelotopes in ${\mathbb R}^n.$ Suppose
$D^*$ is a homothetic copy of $D$
with ratio $\tau>1.$ 
If
$D\subset S \subset D^*,$ 
then $\tau\geq n.$
}
This proposition holds true also for balls.
In fact, the Euler inequality is equivalent to
the following statement.
{\it Suppose $B$ is a ball with radius $r_1$ and
$B^*$ is a ball with radius $r_2$. If 
$B\subset S\subset B^*$, then $r_1\leq nr_2.$ 
Equality takes place if and only if $S$ is a regular simplex
inscribed into $B^*$ and
$B$ is the ball inscribed into $S$.} Another
equivalent form of these propositions is given by
Theorem 2 
(see the note after the proof of this theorem).

Let $x^{(1)},$ $\ldots,$
$x^{(n+1)}$ be the vertices and let $\lambda_1,$ $\ldots,$
$\lambda_{n+1}$ be the basic Lagrange polynomials of an nondegenerate simplex 
$S\subset {\mathbb R}^n$ (see Section 1). 
In what follows $\Gamma_j$ is the $(n-1)$-dimensional hyperplane given
by the equation $\lambda_j(x)=0$, by 
$\Sigma_j$ we mean the $(n-1)$-dimensional face of $S$ contained
in $\Gamma_j$, 
symbol $h_j$ denotes the height of $S$ conducted from the vertex $x^{(j)}$ 
onto~$\Gamma_j$,
and $r$ denotes the inradius of $S$. Define
$\sigma_j$ as $(n-1)$-measure of $\Sigma_j$ and put
$\sigma:=\sum\limits_{j=1}^{n+1} \sigma_j$.
Consider the vector $a_j:=\{l_{1j},\ldots,l_{nj}\}$. This vector is orthogonal
to $\Gamma_j$ and directed into the subspace containing
$x^{(j)}$. Obviously,
$$\lambda_j(x)=
l_{1j}x_1+\ldots+
l_{nj}x_n+l_{n+1,j}=(a_j,x)+l_{n+1,j}=(a_j,x)+\lambda_j(0).$$


\smallskip
{\bf Theorem 1.}
{\it The following equalities are true:
\begin{equation}\label{alpha_bs_sum_l_ij_equality}
\alpha(B_n;S)=
\sum_{j=1}^{n+1}\left(\sum_{i=1}^n l_{ij}^2\right)^{1/2},
\end{equation}
\begin{equation}\label{alpha_bs_h_j_equality}
\alpha(B_n;S)=\sum_{j=1}^{n+1}\frac{1}{h_j},
\end{equation}
\begin{equation}\label{alpha_bs_1_r_equality}
\alpha(B_n;S)= \frac{1}{r},
\end{equation}
\begin{equation}\label{alpha_bs_sigma_nV}
\alpha(B_n;S)=\frac{\sigma}{n\vo(S)}.
\end{equation}
}

\smallskip
{\it Proof.} Let us obtain these pairwise-equivalent equalities from the top up to 
the bottom. 
First we note that 
\begin{equation}\label{alpha_bs_equality}
\alpha(B_n;S)=
\sum_{j=1}^{n+1} \max_{x\in B} (-\lambda_j(x))+1.
\end{equation}
Formula (\ref{alpha_bs_equality})
is the particular case of (\ref{alpha_cs_equality}) in the situation $C=B_n$.
By the Cauchy inequality, 
\begin{equation}\label{cauchy_ineq}
-\|a_j\|\|x\|\leq (a_j,x)\leq
\|a_j\|\|x\|,
\end{equation}
$$-\|a_j\|\|x\|-\lambda_j(0)\leq -\lambda_j(x)\leq
\|a_j\|\|x\|-\lambda_j(0).$$
Both the upper and the lower bounds in 
(\ref{cauchy_ineq}) are reachable. This gives 
$$\max_{x\in B_n} (-\lambda_j(x))=
\max_{\|x\|\leq 1} (-\lambda_j(x))=
\|a_j\|-\lambda_j(0).$$
Therefore,
$$\alpha(B_n;S)=
\sum_{j=1}^{n+1} \max_{x\in B_n} (-\lambda_j(x))+1=
\sum_{j=1}^{n+1}\|a_j\|-\sum_{j=1}^{n+1}\lambda_j(0)+1= 
\sum_{j=1}^{n+1}\left(\sum_{i=1}^n l_{ij}^2\right)^{1/2}.
$$
We made use of the equality 
$\sum\limits_{j=1}^{n+1}\lambda_j(0)=1.$
Since $\lambda_j\left(x^{(j)}\right)=1$, we have
$$h_j={\rm dist}\left(x^{(j)};\Gamma_j\right)=
\frac{\left|\lambda_j\left(x^{(j)}\right)\right|}{\|a_j\|}=
\frac{1}{\|a_j\|}=\frac{1}{\left(\sum\limits_{i=1}^n l_{ij}^2\right)^{1/2}}.$$
Consequently,
$$\alpha(B_n;S)=
\sum_{j=1}^{n+1}\left(\sum_{i=1}^n l_{ij}^2\right)^{1/2}
=\sum_{j=1}^{n+1}\frac{1}{h_j}.$$
We have obtained both (\ref{alpha_bs_sum_l_ij_equality})
and (\ref{alpha_bs_h_j_equality}).

Let us prove
(\ref{alpha_bs_1_r_equality}). 
The ball $B_n$ is a subset of a tranlate
of the simplex $\alpha(B_n;S)S$. This means that a translate of the ball
$\frac{1}{\alpha(B_n;S)}B_n$ 
is contained in $S$. Since the maximum of the radii of balls being contained
in $S$ is equal to $r$, holds true 
$\frac{1}{\alpha(B_n;S)}\leq r,$ i.\,e.,
$\alpha(B_n;S)\geq \frac{1}{r}$. 
To obtaine the inverse inequality, denote by
$B^\prime$ a ball of radius $r$ inscribed into $S$. Then the ball
$B_n=\frac{1}{r}B^\prime$ 
is a subset of some translate of 
$\frac{1}{r}S$.
Using the definition of $\alpha(B_n;S)$ we can write 
$\alpha(B_n;S)\leq \frac{1}{r}$.  
So, we have $\alpha(B_n;S)=\frac{1}{r}$.  

Finally, in order to establish
(\ref{alpha_bs_sigma_nV}), it is sufficient to
utilize (\ref{alpha_bs_1_r_equality}) 
and the formula $\vo(S)=\frac{1}{n}\sigma r$. The latter equality
one can obtain from an ordinary formula for the volume of a simplex
after subdividing $S$ onto $n+1$ simplices in such a way that $j$th of these simplices
has a vertex in the center of the inscribed ball and is supported on  $\Sigma_j$.
\hfill$\Box$

\smallskip
{\bf Corollary 1.}
{\it  We have 
$$\frac{1}{r}=\sum_{j=1}^{n+1}\frac{1}{h_j}.$$
}

\smallskip
For proving, it is sufficient to apply
(\ref{alpha_bs_h_j_equality}) and
(\ref{alpha_bs_1_r_equality}). 
It seems to be interesting
that this geometric relation (which evidently can be obtained also
in a direct way) occurs to be equivalent to general formula for
$\alpha(C;S)$ in the particular 
case when a conveх body $C$ coincide with an Euclidean unit ball.

\smallskip
{\bf Corollary 2.}
{\it  The inradius $r$ and the incenter $z$ of a simplex $S$ can be calculated
by the following formulae:
\begin{equation}\label{r_formula}
r=\frac{1}{ \sum\limits_{j=1}^{n+1}\left(\sum\limits_{i=1}^n l_{ij}^2\right)^{1/2}},
\end{equation}
\begin{equation}\label{z_formula}
z=\frac{1}{ \sum\limits_{j=1}^{n+1}\left(\sum\limits_{i=1}^n l_{ij}^2\right)^{1/2}}
\sum\limits_{j=1}^{n+1}\left(\sum\limits_{i=1}^n l_{ij}^2\right)^{1/2} x^{(j)}.
\end{equation}
The tangent point of the ball $B(z;r)$ and facet 
$\Sigma_k$ has the form
\begin{equation}\label{y_k_formula}
y^{(k)}=\frac{1}{ \sum\limits_{j=1}^{n+1}\left(\sum\limits_{i=1}^n l_{ij}^2\right)^{1/2}}
\left[\sum\limits_{j=1}^{n+1}\left(\sum\limits_{i=1}^n l_{ij}^2\right)^{1/2} x^{(j)}
-\frac{1}{\left(\sum\limits_{i=1}^n l_{ik}^2\right)^{1/2}} \left(l_{1k},\ldots,l_{nk}\right)
\right].
\end{equation}
}

\smallskip
{\it Proof.} Equality
(\ref{r_formula}) follows immediately from
(\ref{alpha_bs_sum_l_ij_equality}) and
(\ref{alpha_bs_1_r_equality}). To obtain
(\ref{z_formula}), let us remark that
$$r=
{\rm dist}(z;\Gamma_j)=
\frac{|\lambda_j(z)|}{\|a_j\|}.$$
Since $z$ lies inside $S$, each barycentric coordinate of this point 
$\lambda_j(z)$ is positive, i.\,e.,
$\lambda_j(z)=r\|a_j\|.$
Consequently,
$$z=\sum_{j=1}^{n+1}\lambda_j(z)x^{(j)}=
r\sum_{j=1}^{n+1} \|a_j\| x^{(j)}.$$
This coincides with (\ref{z_formula}). 
Finally, since
vector $a_k=\{l_{1k},\ldots,l_{nk}\}$
is orthogonal to 
$\Sigma_k$ and is directed
from this facet inside the simplex, a unique common point of
$B(z;r)$ and $\Sigma_k$ has the form
$$y^{(k)}=z-\frac{r}{\|a_k\|}a_k=r\left( \sum_{j=1}^{n+1} \|a_j\| x^{(j)}-\frac{1}{\|a_k\|}
a_k\right).$$
The latter is equivalent to (\ref{y_k_formula}).
\hfill$\Box$

\smallskip
It is interesting to compare 
(\ref{alpha_bs_sum_l_ij_equality}) with the formula 
(\ref{alpha_q_prime_s_formula}) for $\alpha(Q_n^\prime;S)$. Since $B_n$ is
a subset of the cube
$Q_n^\prime=[-1,1]^n$, we have $\alpha(B_n;S)\leq \alpha(Q_n^\prime;S)$.
Analytically, this also follows from the estimate
$$\left(\sum_{i=1}^n l_{ij}^2\right)^{1/2}\leq 
\sum_{i=1}^n |l_{ij}|.$$

For arbitrary $x^{(0)}$ and $\varrho>0$, the number
$\alpha\left(B(x^{(0)};\varrho);S\right)$ can be calculated
with the use of Theorem 1
and the equality 
$\alpha(B(x^{(0)};\varrho);S)$ $=$ $\varrho\alpha(B_n;S)$.

If $S\subset Q_n$, then all the axial diameters
$d_i(S)$ do not exceed $1$ and
(\ref{alpha_d_i_formula}) immediately gives $\alpha(Q_n;S)\geq n$. Moreover,
the equality $\alpha(Q_n;S)=n$ holds when and only when each 
$d_i(S)=1$.
The following proposition
expresses the analogues of these properties for simplices contained
in a ball.  
  
\smallskip
{\bf Theorem 2.}
{\it If $S\subset B_n$, then $\alpha(B_n;S)\geq n.$ The equality
$\alpha(B_n;S)=n$ holds true if and only if 
$S$ is a regular simplex inscribed into $B_n$.
}

\smallskip
{\it Proof.} 
By the definition of $\alpha(B_n;S)$, the ball $B_n$ 
is contained in a translate of the simplex
$\alpha(B_n;S)S$. Hence, some translate 
$B^\prime$  of the ball
$\frac{1}{\alpha(B_n;S)}B_n$  is a subset of
$S$. So, we have the inclusions $B^\prime\subset S\subset B_n$. Since the radius of
$B^\prime$ is equal to
$\frac{1}{\alpha(B_n;S)}$, the inradius
$r$ and the circumradius $R$ of $S$ satisfy the inequalities 
$\frac{1}{\alpha(B_n;S)}\leq r,$ $R\leq 1$. 
Making use of the Euler inequality 
$R\geq nr$, we can write
\begin{equation}\label{theor2_ineqs}
\frac{1}{\alpha(B_n;S)}\leq r\leq \frac{R}{n}\leq \frac{1}{n}. 
\end{equation}
Therefore, $\alpha(B_n;S)\geq n.$ 

The equality $\alpha(B_n;S)=n$ means that the left-hand value in
(\ref{theor2_ineqs}) coincides with the right-hand one. Thus, 
all the inequalities in this chain turn into equalities. We obtain
$R=1,$ $r=\frac{1}{n}$.
Since in this case the Euler inequality 
(\ref{euler_ineq}) 
also becomes an equality, 
$S$ is a regular simplex inscribed into $B_n$.
Conversely, if $S$ is a regular simplex inscribed
into $B_n$, then
$r=\frac{1}{n}$, i.\,e., $\alpha(B_n;S)=\frac{1}{r}=n$.
\hfill$\Box$

\smallskip
We see that Theorem 2 follows from the Euler inequality 
(\ref{euler_ineq}). In fact, these statements are equivalent.
Indeed, suppose $S$ is an arbitrary $n$-dimensional simple,
$r$ is the inradius and $R$ is the circumradius of $S$. 
Let us denote by $B$ the ball containing $S$ and having radius $R$. 
Then some translate $S^\prime$ of the simplex $\frac{1}{R}S$
is contained in $B_n$.
By Theorem~1,
$\alpha(B_n;S^\prime)$ is the inverse to the inradius
of $S^\prime$, i.\,e., is equal to $\frac{R}{r}.$ 
Now assume that Theorem~2 is true.
Let us apply this theorem to the simplex
$S^\prime\subset B_n$. This gives
$\alpha(B_n;S^\prime)=\frac{R}{r}\geq n$
and we have (\ref{euler_ineq}). Finally, if $R=nr,$ then
$\alpha(B_n;S^\prime)=n$. From Theorem 2 we obtain that both
$S^\prime$ and 
 $S$ are regular simplices.

It follows from  
(\ref{ksi_alpha_n_ineq}) that the minimum value of
$\alpha(Q_n;S)$ for $S\subset Q_n$ also is equal to
$n$. 
This minimal value corresponds to those and only those 
$S\subset Q_n$ for which every axial diameter
$d_i(S)$ is equal to $1$. 
The noted property is fulfilled for
the maximum volume simplices in $Q_n$ 
(see [3]), but not for the only these simplices,
if $n>2$.

\section{The value $\xi(B_n;S)$}\label{nev_s3}
In this section, we will obtain the computational
formula for 
the absorption index 
of a simplex $S$ with respect to an Euclidean ball.
We use the previous denotations.

\smallskip
{\bf Theorem 3.}
{\it Suppose $S$ is a nondegenerate simplex in
${\mathbb R}^n$, $x^{(0)}\in {\mathbb R}^n$,
$\varrho>0$. 
If $B\left(x^{(0)};\varrho\right)
 \not\subset S$, we have
\begin{equation}\label{ksi_b_x0_ro_s_l_ij_equality}
\xi\left(B\left(x^{(0)};\varrho\right);S\right)=
(n+1)\max_{1\leq j\leq n+1}
\left[\varrho\left(\sum_{i=1}^n l_{ij}^2\right)^{1/2}-
\sum_{i=1}^n l_{ij}x_i^{(0)}-l_{n+1,j}\right]+1.
\end{equation}
In particular, if
$B_n\not\subset S$, then
\begin{equation}\label{ksi_bs_l_ij_equality}
\xi(B_n;S)=
(n+1)\max_{1\leq j\leq n+1}\left[\left(\sum_{i=1}^n l_{ij}^2\right)^{1/2}
-l_{n+1,j}\right]+1.
\end{equation}
}

\smallskip
{\it Proof.} 
Let us apply the general formula  
(\ref{ksi_cs_equality}) 
in the case $C=B\left(x^{(0)};\varrho\right)$.
The Cauchy inequality yields 
\begin{equation}\label{cauchy_for_ksi_ineq}
-\|a_j\|\|x-x^{(0)}\|\leq (a_j,x-x^{(0)})\leq
\|a_j\|\|x-x^{(0)}\|.
\end{equation}
If $\|x-x^{(0)}\|\leq \varrho$, we see that
$$
-\varrho\|a_j\|\leq 
(a_j,x)-(a_j,x^{(0)})
\leq  \varrho\|a_j\|,
$$
$$-\lambda_j(x)=-(a_j,x)-l_{n+1,j}\leq
\varrho\|a_j\|-(a_j,x^{(0)})-l_{n+1,j}.$$
Since both the upper and the lower bounds  in
(\ref{cauchy_for_ksi_ineq}) are reachable,  
$$\max_{\|x-x^{(0)}\|\leq \varrho} (-\lambda_(x))=
\varrho\left(\sum_{i=1}^n l_{ij}^2\right)^{1/2}-
\sum_{i=1}^n l_{ij}x_i^{(0)}-l_{n+1,j}.$$
It follows that   
$$\xi\left(B\left(x^{(0)};\varrho\right);S\right)=
(n+1)\max_{1\leq j\leq n+1, \|x-x^{(0)}\|\leq \varrho} (-\lambda_j(x))+1=$$
$$=(n+1)\max_{1\leq j\leq n+1}
\left[\varrho\left(\sum_{i=1}^n l_{ij}^2\right)^{1/2}-
\sum_{i=1}^n l_{ij}x_i^{(0)}-l_{n+1,j}\right]+1,$$
and we obtain 
(\ref{ksi_b_x0_ro_s_l_ij_equality}). Equality
(\ref{ksi_bs_l_ij_equality}) appears from 
(\ref{ksi_b_x0_ro_s_l_ij_equality}) for $x^{(0)}=0, \varrho=1$.
\hfill$\Box$

\section{The equality $\beta_n=n$. Commentaries}\label{nev_s4}

\smallskip
{\bf Theorem 4.}
{\it 
If $S\subset B_n$, then $\xi(B_n;S)\geq n.$ The equality
$\xi(B_n;S)=n$ takes place if and only if 
$S$ is a regular simplex inscribed into $B_n$.
}

\smallskip
{\it Proof.} The statement immediately follows from Theorem 2 
and the inequality
$\xi(B_n;S)\geq \alpha(B_n;S)$. We give
here also a direct proof without applying the Euler inequality
that was used to obtain the estimate 
$\alpha(B_n;S)\geq n$.  

First let $S$ be a regular simplex inscribed into $ B_n $. 
Then $\alpha(B_n;S)=n$ and the inradius of
$S$ is equal to $\frac{1}{n}$. 
Since the simplex $\xi(B_n;S)S$ is circumscribed around $B_n$, we have
the equalities
$\xi(S;B_n)=\alpha(S;B_n)=n$ and also
relation (\ref{relation_cs}) with $C=B_n$.
It follows from (\ref{ksi_cs_equality}) that for any 
$j=1,\ldots,n+1$  
$$\max_{x\in B_n} (-\lambda_j(x))=\frac{n-1}{n+1},$$
where  $\lambda_j$ are the basic Lagrange polynomials
related to $S$.

Now suppose simplex $S$ is contained in $B_n$ but is not regular or is not
inscribed into the ball.
Denote the Lagrange polynomials of this simplex by
$\mu_j$. There exist a regular simplex 
$S^*$ inscribed into $B_n$ and an integer $k$ such that
$S$ is contained in the strip $0\leq\lambda_k(x)\leq 1$, the $k$th $(n-1)$-dimensional
faces of $S$ and $S^*$ are parallel, and $S$ has not any common points with
at least one of the boundary hyperplanes of this strip.
Here $\lambda_j$ are the basic Lagrange polynomials of $S^*$. 
The vertex $x^{(k)}$ of the simplex $S^*$ does not lie in its
$k$th facet. 
Assume $u$ is a point of the boundary of $B_n$
most distant from $x^{(k)}$. Then
$u$ is the maximum point of polynomial $-\lambda_k(x)$, i.\,e.,
 $- \lambda_j(u)=\frac{n-1}{n+1}$. 
Consider the straight line connecting 
$x^{(k)}$ and $u$. 
Denote by $y,z$ and $t$ the inersection points 
of this line and pairwize parallel hyperplanes
$\mu_k(x)=1,$ $\mu_k=0$ and $\lambda_k(x)=0$ respectively.
We have 
\begin{equation}\label{ineqs_one_strong}
\|x^{(k)}-t\|\geq \|y-z\|, \quad \|t-u\|\leq \|z-u\|.
\end{equation} 
At least one of these inequalities is fulfilled
in the strict form.
The linearity of the basic Lagrange polynomials means that
$$\frac{\mu_k(z)-\mu_k(u)}{\mu_k(y)-\mu_k(z)}=
\frac{\|z-u\|}{\|y-z\|}, \quad 
\frac{\lambda_k(t)-\lambda_k(u)}{\lambda_k\left(x^{(k)}\right)-\lambda_k(t)}=
\frac{\|t-u\|}{\left\|x^{(k)}-t\right\|}.$$
Since $\mu_k(y)=1,$ $\mu_k(z)=0,$ $\lambda_k\left(x^{(k)}\right)=1,$ and $\lambda_k(t)=0$, 
we get
$$-\mu_k(u)=  \frac{\|z-u\|}{\|y-z\|} > 
\frac{\|t-u\|}{\left\|x^{(k)}-t\right\|}=-\lambda_k(u)=\frac{n-1}{n+1}.$$
We made use of (\ref{ineqs_one_strong}) and took into account that
at least one of the inequalities is strict.
The application of (\ref{ksi_cs_equality}) yields
$$\xi(B_n;S)=(n+1)\max_{1\leq j\leq n+1}
\max_{x\in B_n}(-\mu_j(x))+1\geq (n+1)(-\mu_k(u))+1>n.$$ 
Thus, if $S$ is not regular simplex inscribed into
$B_n$, then
$\xi(B_n;S)>n$. 

We see that each simplex $S\subset B_n$ satisfies the estimate 
$\xi(B_n;S)\geq n$. The equality takes place if and only if
$S$ is a regular simplex inscribed into $B_n$. 
\hfill$\Box$

By analogy with the value 
$\xi_n=\min\{\xi(Q_n;S): S\subset Q_n\}$
defined through the unit cube, let us 
introduce the similar numerical characteristic given by the unit ball:
$$\beta_n:=\min \{ \xi(B_n;S): \,
S \mbox{ --- $n$-мерный симплекс,} \,
S\subset B_n, \, \vo(S)\ne 0\}.$$
Many problems concerning $\xi_n$ yet have not been solved. 
For example, $\xi_2 = 1 + \frac{3\sqrt{5}}{5}$ still remains
the only accurate value of 
$\xi_n $ for even $n$; moreover,
this value was discovered in a rather difficult way
(see [3; Chapter\,2]).
Compared to $\xi_n$
the problem on numbers $\beta_n $ turns out to be trivial.

\smallskip
{\bf Corollary 3.}
{\it  For any $n$, we have $\beta_n=n$. 
The only simplex $S\subset B_n$  extremal with respect to $\beta_n$   
is an arbitrary regular simplex inscribed into $B_n$.}

\smallskip
{\it Proof.} It is sufficient to apply Theorem 4. \hfill$\Box$

\smallskip 
The technique developed for a ball makes it possible
to illustrate some results 
having been earlier got for a cube.
Here we note  a proof of the following known statement  
which differs from the proofs given in
[3; \S\,3.2] and [12].

\smallskip
{\bf Corollary 4.}
{\it If
there exists an Hadamard matrix
of order $n+1$, then $\xi_n=n.$
}

\smallskip
{\it Proof.} It is known (see, e.\,g., [9]) 
that for these and only these~$n$ we can inscribe into $Q_n$ a regular simplex $S$
so that all the vertices of $S$ will coincide with vertices of
the cube. Let us denote by $B$ the ball with radius $\frac{\sqrt{n}}{2}$
having the center in center
of the cube.
Clearly, $Q_n$ is inscribed into $B$, therefore,
the simplex is inscribed into the ball as well.
Since $S$ is regular, by Theorem 4 and by similarity reasons,
we have
$\xi(B;S)=n.$ 
The inclusion $Q_n\subset B$ means that
$\xi(Q_n;S)\leq \xi(B;S),$ i.\,e.~$\xi(Q_n;S)\leq n$.
From (\ref{ksi_alpha_n_ineq}) it follows that the inverse inequality 
$\xi(Q_n;S)\geq n$ is also true. Hence,
$\xi(Q_n;S)=n$. Simultaneously 
(\ref{ksi_alpha_n_ineq}) gives
$\xi_n=\xi(Q_n;S)=n$. 
\hfill$\Box$

This argument is based on the following fact:
if $S$ is a regular simplex with the vertices in vertices of
$Q_n$, then the simplex $nS$ absorbs not only the cube $Q_n$ but also
the ball $B$ circumscribed  around the cube.
The corresponding absorption index $n$ is the minimum possible
both for the cube and the ball. In addition, we mention the following
property.

\smallskip
{\bf Corollary 5.}
{\it Assume that $S\subset Q_n\subset nS$
and simplex $S$ is not regular. Then
$B\not\subset nS$. 
}

\smallskip
{\it Proof.}   The inclusion $B\subset nS$ implies that   $\xi(B;S)=n$.  This
way $S$ is a regular simplex inscribed into the ball
$B$. But since this is not so,  
$B$ is not a subset of $nS$.
\hfill$\Box$

\smallskip
Simplices satisfying the condition of Corollary 5 
exist at least for
$n=3, 5,$ and $9$ (see [12]).

The relations (\ref{ksi_alpha_n_ineq}) mean that always
$\xi_n\geq n$. Since $\xi_2=1+\frac{3\sqrt{5}}{5}>2$,
there exist $n$'s such that
$\xi_n>n$. 
Besides the cases when $n+1$ is an Hadamard number,
the equality $\xi_n=n$ is established for $n=5$ and $n=9$ 
(the extremal simplices in
${\mathbb R}^5$ and ${\mathbb R}^9$ are given in
[12]).
For all such dimensions holds true  
$\xi_n=\beta_n$, i.\,e., with respect to the minimum
absorption index of an internal simplex,
both the convex bodies, an $n$-dimensional cube and
an $n$-dimensional ball, have the same behavoir.

The equality $\xi_n=n$ is equivalent to the existence
of simplices satisfying the inclusions $S\subset Q_n\subset nS$.
Some properties of such simplices (e.\,g., the fact that the center of gravity of
$S$ coincides with the center of the cube; see [7])
are similar to the properties of regular simplices
inscribed into the ball.
However,
the problem to describe the set of all dimensions 
where exist those simplices, seems to be very difficult
and nowaday is far from  solution.

  \bigskip
\centerline{\bf\Large References}

\begin{itemize}
\item[1.]
 Nevskij,~M.\,V., On a certain relation for the minimal norm of an interpolational 
projection,
{\it Model.  Anal. Inform. Sist.}, 2009, vol.~16, no.~1, pp.~24--43 (in~Russian).

\item[2.]
 Nevskii,~M.\,V. On a property of $n$-dimensional simplices,
{\it Math. Notes}, 2010, vol.~87, no.~4, pp.~543--555.

\item[3.]
 Nevskii,~M.\,V.,
{\it Geometricheskie ocenki v polinomialnoi interpolyacii} 
(Geometric Estimates in Polynomial
Interpolation), Yaroslavl': Yarosl. Gos. Univ., 2012 (in~Russian).
.

\item[4.]
Nevskii,~M.\,V., On the minimal positive homothetic image of a simplex containing a 
convex body,
{\it Math. Notes}, 2013, vol.~93, no.~3--4, pp.~470--478.

\item[5.]
Nevskii,~M.\,V., and Ukhalov, A.\,Yu.,
On numerical charasteristics of a simplex and their estimates,
{\it Model. Anal. Inform. Sist.}, 2016, vol.~23, no.~5, pp.~603--619 
(in~Russian).
English transl.: {\it Aut.
Control Comp. Sci.}, 2017, vol.~51, no.~7, pp.~757--769.

\item[6.]
Nevskii,~M.\,V., and Ukhalov, A.\,Yu.,
New estimates of numerical values related to a simplex,
{\it Model. Anal. Inform. Sist.}, 2017, vol.~24, no.~1, pp.~94--110
(in~Russian).
English transl.: {\it Aut.
Control Comp. Sci.}, 2017, vol.~51, no.~7, pp.~770--782.

\item[7.]
Nevskii,~M.\,V., and Ukhalov, A.\,Yu.,
On $n$-dimensional simplices satisfying inclusions 
$S\subset [0,1]^n\subset nS$,
{\it Model. Anal. Inform. Sist.}, 2017, vol.~24, no.~5, pp.~578--595
(in~Russian).
English transl.: {\it Aut.
Control Comp. Sci.}, 2018, vol.~52, no.~7, pp.~667--679.

\item[8.]
Nevskii,~M.\,V., and Ukhalov, A.\,Yu.,
On minimal  absorption index for an $n$-dimensional simplex,
{\it Model. Anal. Inform. Sist.}, 2018, vol.~25, no.~1, pp.~140--150
(in~Russian).
English transl.: {\it Aut.
Control Comp. Sci.}, 2018, vol.~52, no.~7, pp.~680--687.

\item[9.]
Hudelson,~M., Klee, V., and Larman,~D.,
Largest $j$-simplices in $d$-cubes: some relatives of the
Hadamard maximum determinant problem,
{\it Linear Algebra Appl.}, 1996, vol.~241--243, pp.~519--598

\item[10.]
Klamkin~M.\,S., and Tsifinis~G.\,A.,
Circumradius--inradius inequality for a simplex,
{\it Mathematics Magazine}, 1979, vol.~52, no.~1, pp.~20--22.

\item[11.]
Nevskii,~M.,
Properties of axial diameters of a simplex,
{\it Discr. Comput. Geom.}, 2011, vol.~46, no.~2, pp.~301--312.

\item[12.]
Nevskii,~M., and Ukhalov A.,
Perfect simplices in ${\mathbb R}^5$,
{\it 
Beitr. Algebra Geom.},
2018, vol.~59, no.~3, pp.~501--521.

\item[13.]
Yang~S., and Wang~J.,
Improvements of $n$-dimensional Euler inequality,
{\it Journal of~Geometry}, 1994, vol.~51, pp.~190--195

\item[14.]
Vince~A., A simplex contained in a sphere,
{\it Journal of~Geometry}, 2008, vol.~89, no.~1--2, 
pp.~169--178.

\end{itemize}

\end{document}